\documentclass[11pt]{amsart}
\usepackage{amssymb}
\usepackage{amsmath}
\usepackage[dvips]{graphicx}
\usepackage{psfrag}
\begin{document}
\setlength{\parindent}{1.2em}
\def\COMMENT#1{}
\def\TASK#1{}
\def\noproof{{\unskip\nobreak\hfill\penalty50\hskip2em\hbox{}\nobreak\hfill%
       $\square$\parfillskip=0pt\finalhyphendemerits=0\par}\goodbreak}
\def\endproof{\noproof\bigskip}
\newdimen\margin   
\def\textno#1&#2\par{%
   \margin=\hsize
   \advance\margin by -4\parindent
          \setbox1=\hbox{\sl#1}%
   \ifdim\wd1 < \margin
      $$\box1\eqno#2$$%
   \else
      \bigbreak
      \hbox to \hsize{\indent$\vcenter{\advance\hsize by -3\parindent
      \sl\noindent#1}\hfil#2$}%
      \bigbreak
   \fi}
\def\proof{\removelastskip\penalty55\medskip\noindent{\bf Proof. }}
\def\enddiscard{}
\long\def\discard#1\enddiscard{}
\newtheorem{firstthm}{Proposition}
\newtheorem{thm}[firstthm]{Theorem}
\newtheorem{prop}[firstthm]{Proposition}
\newtheorem{lemma}[firstthm]{Lemma}
\newtheorem{cor}[firstthm]{Corollary}
\newtheorem{problem}[firstthm]{Problem}
\newtheorem{defin}[firstthm]{Definition}
\newtheorem{conj}[firstthm]{Conjecture}
\newtheorem{theorem}[firstthm]{Theorem}
\newtheorem{claim}{Claim}
\def\eps{{\varepsilon}}
\def\N{\mathbb{N}}
\def\R{\mathbb{R}}
\def\P{\mathcal{P}}
\def\Ybar{\bar{Y}}

\title{Linkedness and ordered cycles in digraphs}
\author{Daniela K\"uhn \and Deryk Osthus}
\date{}
\maketitle \vspace{-.8cm}
\begin{abstract} \noindent
Given a digraph~$D$, let $\delta(D):=\min\{\delta^+(D), \delta^-(D)\}$ be the
minimum degree of~$D$. We show that every sufficiently large digraph~$D$ with $\delta(D)\ge
n/2 +\ell-1$ is $\ell$-linked. The bound on the minimum degree is
best possible and confirms a conjecture of
Manoussakis~\cite{Manoussakis}. We also determine the smallest minimum degree
which ensures that a sufficiently large digraph~$D$ is $k$-ordered, i.e.~that
for every sequence $s_1,\dots,s_k$ of distinct vertices of~$D$ there is a directed cycle
which encounters $s_1,\dots,s_k$ in this order. 
\end{abstract}

\section{Introduction}\label{intro}

The \emph{minimum degree~$\delta(D)$} of a digraph~$D$ is the minimum of its
minimum outdegree~$\delta^+(D)$ and its minimum indegree~$\delta^-(D)$.
When referring to paths and cycles in digraphs we always mean that these are directed without
mentioning this explicitly. 
A digraph $D$ is $\ell$-linked if $|D|\ge 2\ell$ and if for every sequence
$x_1,\dots,x_\ell,y_1,\dots,y_\ell$ of distinct vertices there are disjoint paths
$P_1,\dots,P_\ell$ in~$D$ such that $P_i$ joins~$x_i$ to~$y_i$. Since this is a very strong and 
useful property to have in a digraph, the question of course arises how it can be forced by
other properties.

In the case of (undirected) graphs, much progress has been made in this direction.
In particular, linkedness is closely related to connectivity:
Bollob\'as and Thomason~\cite{BT96} showed that every $22k$-connected graph is $k$-linked
(this was recently improved to~$10k$ by Thomas and Wollan~\cite{ThomasWollan}).
However, for digraphs the situation is quite different:
Thomassen~\cite{Thomassen} showed that for all $k$ there are strongly $k$-connected digraphs which
are not even $2$-linked.


Our first result determines the minimum degree forcing a (large) digraph to be $\ell$-linked,
which confirms a conjecture of Manoussakis~\cite{Manoussakis} for large digraphs. 

\begin{thm}\label{thmlinked}
Let $\ell\ge 2$.
Every digraph $D$ of order $n\ge 1600\ell^3$ which satisfies $\delta(D)\ge n/2 +\ell-1$
is $\ell$-linked.
\end{thm}
It is not hard to see that the bound on minimum degree in Theorem~\ref{thmlinked} is best possible
(see Proposition~\ref{propbestlinked}). It is also easy to see that for $\ell=1$ the correct bound
is $\delta(D)\ge \lfloor n/2 \rfloor$.
The cases~$\ell=2,3$ of Theorem~\ref{thmlinked}
were proved by Heydemann and Sotteau~\cite{HS} and Manoussakis~\cite{Manoussakis}
respectively. Manoussakis~\cite{Manoussakis} also determined the number of edges which force a digraph to be
$\ell$-linked. A discussion of these and related results can be found in the monograph 
by Bang-Jensen and Gutin~\cite{digraphsbook}.

Note that it does not make sense to ask for the
minimum outdegree of a digraph~$D$ which ensures that~$D$ is $\ell$-linked
(or similarly, to ask for the minimum indegree). Indeed, the digraph obtained from a
complete digraph~$A$ of order~$n-1$ by adding a new vertex~$x$ which sends an edge
to every vertex in~$A$ has minimum outdegree $n-2$ but is not even 1-linked.

A slightly weaker notion is that of a $k$-ordered digraph:
a digraph $D$ is \emph{$k$-ordered} if $|D|\ge k$ and if for every sequence $s_1,\dots,s_k$
of distinct vertices of~$D$ there is a cycle which encounters $s_1,\dots,s_k$ in this order.
It is not hard to see that every $\ell$-linked digraph is also $\ell$-ordered.%
     \COMMENT{This is true if every $\ell$-linked digraph $D$ has minimum outdegree at least~$2\ell-1$.
(Indeed, if the latter is true then for each $s_i$ we can pick an outneighbour $s'_i$ such that
$s_1,\dots,s_\ell,s'_1,\dots,s'_\ell$ are distinct. If $D$ is $\ell$-linked then we can find disjoint
$s'_i$-$s_{i+1}$ paths which together with the edges $\overrightarrow{s_is'_i}$ form the required
cycle.) If some vertex $x$ has outdegree $\le 2\ell-2$ then $D$ is not $\ell$-linked:
take $x_1:=x$, let $x_2,\dots,x_\ell$ and $y_2,\dots,y_\ell$ be the $2\ell-2$ vertices in the
outneighbourhood of~$x$ and let $y_1$ be any other vertex.}
Conversely, every~$2\ell$-ordered digraph~$D$ is also $\ell$-linked: if $x_1,\dots,x_\ell,y_1,\dots,y_\ell$
is a sequence of vertices as in the definition of $\ell$-linkedness then a cycle which encounters
$x_1,y_1,x_2,y_2,\dots,x_\ell,y_\ell$ in this order would yield the paths required
for the linking. The next result says that as far as the minimum degree is concerned
it is no harder to
guarantee the~$2\ell$ paths forming such a cycle than to guarantee just the~$\ell$ paths
required for the linking. In particular, note that Theorem~\ref{thmkordered} immediately
implies Theorem~\ref{thmlinked}.
\begin{thm}\label{thmkordered}
Let $k\ge 2$.
Every digraph $D$ of order $n\ge 200k^3$ which satisfies $\delta(D)\ge  (n+k)/2-1$
is $k$-ordered. 
\end{thm}
Again, the bound on the minimum degree is best possible (see Proposition~\ref{propbestordered}).
Moreover, it is easy to see that if $k=1$ then the correct bound is $\delta(D)\ge n/2 -1$.
The proof of Theorem~\ref{thmkordered} yields paths between the $k$ `special' vertices whose 
length is at most~6 and it is also easy to translate the proof into an algorithm which finds these
paths in polynomial time (see the remarks after the end of the proof).

Somewhat surprisingly, the minimum degree in both theorems is not quite the same as in
the undirected case: Kawarabayashi, Kostochka and Yu~\cite{KKYlinked} proved that
the smallest minimum degree which guarantees a 
graph on $n$ vertices to be $\ell$-linked is $\lfloor n/2 \rfloor +\ell-1$ for large $n$.
(Egawa et al.~\cite{Egawaetal} independently determined the smallest minimum degree
which guarantees the existence of $\ell$ disjoint cycles containing $\ell$ specified independent 
edges, which is clearly a very similar property.)
Kierstead, Sark\"ozy and Selkow~\cite{KSSkordered} proved that
the smallest minimum degree which guarantees a 
graph on $n$ vertices to be $k$-ordered is 
$\delta(D)\ge \lceil n/2\rceil + \lfloor k/2 \rfloor -1$ for large $n$. So in the
undirected case the `$2\ell$-ordered' result does not quite imply the `$\ell$-linked' result.
The proofs in~\cite{Egawaetal,KKYlinked,KSSkordered} do not seem to generalize to digraphs.

\section{Further work and open problems} \label{openprobs}

In a sequel to this paper, we hope to apply Theorem~\ref{thmkordered} to obtain the 
following stronger results, which would also generalize the theorem of 
Ghouila-Houri~\cite{GhouilaHouri}
that any digraph $D$ on $n$ vertices with $\delta(D)\ge n/2$ contains
a Hamilton cycle:  
we aim to apply Theorem~\ref{thmkordered} to show that if $k\ge 2$ and~$D$ is
a sufficiently large digraph whose minimum degree is as in Theorem~\ref{thmkordered}
then~$D$ is even $k$-ordered Hamiltonian, i.e.~for every sequence $s_1,\dots,s_k$
of distinct vertices of~$D$ there is a Hamilton cycle which encounters $s_1,\dots,s_k$
in this order.
One can use this to prove that the minimum degree condition in Theorem~\ref{thmlinked}
already implies that the digraph~$D$ is Hamiltonian $\ell$-linked,
i.e.~the paths linking the pairs of vertices span the entire vertex set of~$D$.
Note that this in turn would immediately imply that $D$ is $\ell$-arc ordered Hamiltonian,
i.e.~$D$ has a Hamilton cycle which contains any $\ell$ disjoint edges in a given order.
Note that in each case the examples in Section~\ref{extremal} show that the minimum degree condition 
would be best possible.%
     \COMMENT{They work for $\ell$-arc orderedness too. Indeed, if $n$ is even we take $\ell-2$ edges
$\vec{e}_2,\dots,\vec{e}_{\ell-1}$ inside $A\cap B$, the first edge $\vec{e}_1$ joins a vertex
$a_1\in A\setminus B$ to some vertex $x_1\in A\cap B$ and the last edge $\vec{e}_\ell$ will
join the remaining vertex $x_\ell\in A\cap B$ to some vertex $b_\ell \in B\setminus A$.
Then there is no cycle which contains $\vec{e}_1,\dots,\vec{e}_{\ell}$ in that order.
If $n$ is odd we take $\ell-3$ edges $\vec{e}_3,\dots,\vec{e}_{\ell-1}$ inside $X$.
Let $x,x',x''$ be the remaining vertices in~$X$. Let $\vec{e}_1:= \overrightarrow{xx_1}$,
$\vec{e}_2:= \overrightarrow{y_1x_2}$ and $\vec{e}_\ell:= \overrightarrow{y_2x'}$.
Then any cycle which encounters $\vec{e}_1,\dots,\vec{e}_{\ell}$ in that order
must meet $x''$ when going from $x_1$ to $y_1$ and when going from $x_2$ to $y_2$,
which is impossible.}
Undirected versions of these statements were first obtained 
by Kierstead, Sark\"ozy and Selkow~\cite{KSSkordered} and Egawa et al.~\cite{Egawaetal}
respectively (and a common generalization of these in~\cite{Chenetalk-arc}).

For graphs, the concepts `$\ell$-linked' and `$k$-ordered' were generalized to `$H$-linked' by 
Jung~\cite{Jung}: a graph $G$ is \emph{$H$-linked} 
if $G$ contains a subdivision of $H$ with prescribed branch vertices
(so $G$ is $k$-ordered if and only if it is $C_k$-linked). 
The minimum degree which forces a graph to be $H$-linked for an arbitrary $H$ was determined 
in~\cite{FGTW, KYlinked, KYlinked2, GKYlinked}.
Clearly, one can ask similar questions also for digraphs.

Finally, we believe that the bound on $n$ which we require in Theorem~\ref{thmkordered}
(and thus in Theorem~\ref{thmlinked}) can be reduced to one which is linear in $k$.
 
\section{Notation and extremal examples} \label{extremal}

Before we discuss the examples showing that the bounds on the minimum degree in
Theorems~\ref{thmlinked} and~\ref{thmkordered} are best possible, we will introduce
the basic notation used throughout the paper.
A digraph~$D$ is \emph{complete} if every pair of vertices of~$D$ is joined by edges
in both directions. The order~$|D|$ of a digraph~$D$ is the number of its vertices.
We write~$N^+(x)$ for the outneighbourhood
of a vertex~$x$ and $d^+(x):=|N^+(x)|$ for its outdegree. Similarly, we write~$N^-(x)$
for the inneighbourhood of a vertex~$x$ and $d^-(x):=|N^-(x)|$ for its indegree.
We set $d(x):=\min\{d^+(x),d^-(x)\}$.
Given a set~$A$ of vertices of~$D$, we write $N^+_A(x)$ for the set of all outneighbours
of~$x$ in~$A$. $N^-_A(x)$, $d^+_A(x)$ and $d^-_A(x)$ are defined similarly.
Given two vertices $x,y$ of a digraph~$D$,
an \emph{$x$-$y$ path in~$D$} is a directed path which joins~$x$ to~$y$.
Given two disjoint vertex sets~$A$ and~$B$ of~$D$, an~$A$-$B$ edge is an edge
$\overrightarrow{ab}$ where $a\in A$ and $b\in B$.

The following proposition shows that the bound on the minimum degree in Theorem~\ref{thmlinked}
cannot be reduced.

\begin{prop}\label{propbestlinked}
For every $\ell\ge 2$ and every $n\ge 2\ell$ there exists a digraph~$D$ on~$n$ vertices
with minimum degree~$\lceil n/2\rceil +\ell-2$ which is not $\ell$-linked.
\end{prop}
\proof
We will distinguish the following cases.%
     \COMMENT{The correct threshold for~$\ell=1$ is~$\lfloor n/2\rfloor$.
Indeed, suppose that~$D$ is a digraph of minimum degree at least $\lfloor n/2\rfloor$
and that we wish to find an $x$-$y$ path. Clearly we may assume that $\overrightarrow{xy}$
is not an edge and that $N^+(x)\cap N^-(y)=\emptyset$. But this is impossible
since then $n=|D|\ge |\{x,y\}|+|N^+(x)|+ |N^-(y)|\ge 2+2\lfloor n/2\rfloor>n$.
A digraph of minimum degree $\lfloor n/2\rfloor-1$ can be obtained from disjoint
complete graphs of order $\lfloor n/2\rfloor$ and $\lceil n/2\rceil$.}

\smallskip

\noindent
\textbf{Case 1.} \emph{$n$ is even.}

\medskip

\noindent
Let~$D$ be the digraph which consists of complete digraphs~$A$ and~$B$ of order~$n/2 +\ell-1$
which have precisely $2\ell-2$ vertices in common. To see that~$D$ is not $\ell$-linked
let $x_1,\dots,x_{\ell-1},y_1,\dots,y_{\ell-1}$
denote the vertices in~$A\cap B$. Pick some vertex $x_\ell\in A\setminus B$ and some vertex
$y_\ell\in B\setminus A$. Then~$D$ does not contain disjoint paths between~$x_i$ and~$y_i$ for
all~$i=1,\dots,\ell$. The minimum degree of~$D$ is attained by the vertices in
$(A\setminus B)\cup (B\setminus A)$ and thus is as desired.

\medskip

\noindent
\textbf{Case 2.} \emph{$n$ is odd.}

\smallskip

\noindent
In this case, we define~$D$ as follows.
Let~$A$ and~$B$ be disjoint complete digraphs of order $\lceil n/2\rceil-\ell-1$.
Add a complete digraph~$X$ of order~$2\ell-3$ and join all vertices in~$X$ to all vertices
in~$A\cup B$ with edges in both directions.
Add a set $S:=\{x_1,x_2,y_1,y_2\}$ of~4 new vertices such that each vertex in~$S$
is joined to each vertex in~$X$ with edges in both directions.
Moreover, we add all the edges between different vertices in~$S$ except
for~$\overrightarrow{x_1y_1}$ and~$\overrightarrow{x_2y_2}$.
Finally, we connect the vertices in~$S$ to the vertices
in~$A\cup B$ as follows. Both~$x_1$ and~$y_1$ receive edges from every vertex in~$B$ and send
edges to every vertex in~$A$. Additionally, $x_1$ will receive an edge from every vertex in~$A$
and $y_1$ will send an edge to every vertex in~$B$.
Both~$x_2$ and~$y_2$ receive edges from every vertex in~$A$ and send
edges to every vertex in~$B$. Additionally, $x_2$ will receive an edge from every vertex in~$B$
and $y_2$ will send an edge to every vertex in~$A$ (see Figure~1).
\begin{figure}
\centering
\psfrag{1}[][]{$A$}
\psfrag{2}[][]{$B$}
\psfrag{3}[][]{$X$}
\psfrag{4}[][]{$x_1$}
\psfrag{5}[][]{$y_1$}
\psfrag{6}[][]{$x_2$}
\psfrag{7}[][]{$\ y_2$}
\includegraphics[scale=0.45]{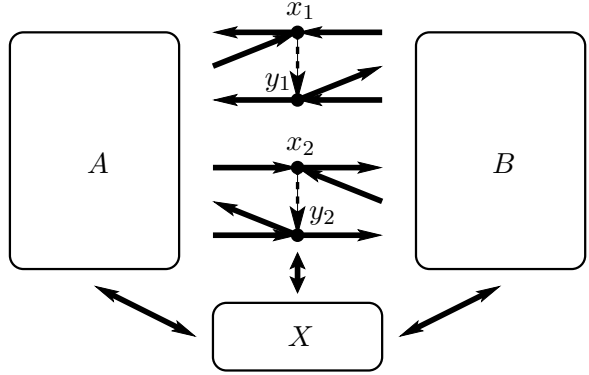}
\caption{The digraph~$D$ in Case~2 of Proposition~\ref{propbestlinked}. 
The dashed arrows indicate the missing edges
between~$x_1$ and~$y_1$ and between~$x_2$ and~$y_2$.}
\end{figure}%

To check that~$D$ has the required minimum degree, consider first any vertex~$a\in A$.
As~$a$ sends edges to~3 vertices in~$S$ and receives edges from~$3$ such vertices, we have
that $d(a)=|A|-1+|X|+3=\lceil n/2\rceil+\ell-2$. It follows similarly that the vertices
in~$B$ have the correct degree. Thus consider any vertex~$s\in S$.
Then~$s$ sends edges to all vertices in~$A$ or to all vertices in~$B$ (or both) and~$s$
receives edges from all vertices in~$A$ or from all vertices in~$B$ (or both).
Thus $d(s)=|A|+|X|+2=\lceil n/2\rceil+\ell-2$. It is easy to check that the vertices in~$X$
have the required degree and thus $\delta(D)=\lceil n/2\rceil+\ell-2$.

To see that~$D$ is not $\ell$-linked, let $x,x_3,\dots,x_\ell,y_3,\dots,y_\ell$
denote the vertices in~$X$. Then we cannot link~$x_i$ to~$y_i$ for each~$i=1,\dots,\ell$
since every $x_1$-$y_1$ path must meet~$X\cup\{x_2,y_2\}$ (and thus would contain~$x$) and the
analogue is true for every~$x_2$-$y_2$ path.
\endproof

We conclude this section with the examples showing that the bound on the minimum
degree in Theorem~\ref{thmkordered} is best possible.

\begin{prop}\label{propbestordered}
For every $k\ge 2$ and every $n\ge 2k$ there exists a digraph~$D$ on~$n$ vertices
with minimum degree~$\lceil (n+k)/2\rceil-2$ which is not $k$-ordered.
\end{prop}
\proof
We will distinguish the following cases.%
      \COMMENT{To see that for every $n\ge 3$ there exists a digraph~$D$ on~$n$ vertices
with minimum degree~$\lceil n/2\rceil-2$ which is not $1$-ordered, look at
the digraph~$D$ obtained from two disjoint complete digraphs~$A$ and~$B$
of order $\lceil n/2\rceil-1$ and $\lfloor n/2 \rfloor$ by adding a new vertex~$s$
and all edges from~$A$ to~$s$ as well as all edges from~$s$ to~$B$.
Then $\delta(D)=|A|-1=\lceil n/2\rceil-2$. $D$ is not~$1$-ordered since~$s$ does
not lie on a cycle.}

\smallskip

\noindent
\textbf{Case 1.} \emph{$k\ge 3$ is odd and~$n$ is even.}

\smallskip

\noindent
In this case, we define~$D$ as follows.
Let~$A$ and~$B$ be disjoint complete digraphs of order $n/2-k+1$.
Add a complete digraph~$X$ of order~$k-2$ and join
all its vertices to all vertices in~$A\cup B$ with edges in both directions.
Add new vertices $s_1,\dots,s_k$ such that every~$s_i$ is joined to all vertices
in~$X$ with edges in both directions. Moreover, we add all the edges $\overrightarrow{s_is_j}$
for~$j\neq i,i+1$ where $s_{k+1}:=s_1$. We also%
     \COMMENT{Adding $\overrightarrow{s_1s_2}$ is not necessary but the picture is
easier to draw that way.}
add the edge $\overrightarrow{s_1s_2}$.
Finally, we connect the~$s_i$ to the vertices
in~$A\cup B$ as follows. Both~$s_1$ and~$s_2$ receive edges from every vertex in~$B$ and send
edges to every vertex in~$B$. Additionally, $s_1$ will send an edge to every vertex in~$A$
and $s_2$ will receive an edge from every vertex in~$A$.
Each of $s_3,s_5,\dots,s_k$ receives an edge from every vertex in~$A$ and
sends an edge to every vertex in~$A$. Each of $s_4,s_6,\dots,s_{k-1}$ receives an edge from
every vertex in~$B$ and sends an edge to every vertex in~$B$ (see Figure~2).
\begin{figure}
\centering
\psfrag{1}[][]{$A$}
\psfrag{2}[][]{$B$}
\psfrag{3}[][]{$X$}
\psfrag{4}[][]{$s_1$}
\psfrag{5}[][]{$s_2$}
\psfrag{6}[][]{$s_3$}
\psfrag{7}[][]{$\ \ s_4$}
\psfrag{8}[][]{$s_5$}
\includegraphics[scale=0.5]{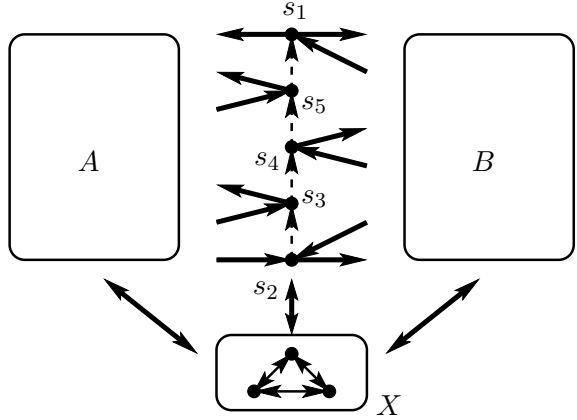}
\caption{The digraph~$D$ for $k=5$ in Case~1 of Proposition~\ref{propbestordered}. 
The dashed arrows indicate missing edges
between the~vertices~$s_i$.}
\end{figure}%
Let us now check that the minimum degree of the digraph~$D$ thus obtained is as required.
Let~$S:=\{s_1,\dots,s_k\}$. Note that each vertex $v\in A\cup B$ sends edges to
precisely $(k+1)/2$ vertices in~$S$ and receives edges from precisely that many vertices.
Since $|A|=|B|$, it follows that $d(v)=|A|-1+|X|+(k+1)/2=n/2-2+(k+1)/2=\lceil (n+k)/2\rceil-2$.
Now consider any~$s_i\in S$. Then~$s_i$ receives edges from either all vertices in~$A$
or all vertices in~$B$ (or both) and $s_i$ sends edges to either all vertices in~$A$
or all vertices in~$B$ (or both). Hence $d(s_i)\ge |A|+|X|+|S|-2=n/2-1+k-2\ge \lceil (n+k)/2\rceil-2$.
It is easy to check that the degree of the vertices in~$X$ is $> \lceil (n+k)/2\rceil-2$.

To see that~$D$ is not~$k$-ordered note that every cycle in~$D$ which
encounters $s_1,\dots,s_k$ in this order would use at least one vertex from~$X$ between~$s_i$
and~$s_{i+1}$ for every $i\neq 1$ (see Figure~2). But since $|X|=k-2$ this is impossible.

\medskip

\noindent
\textbf{Case 2.} \emph{$k$ is even.}

\smallskip

\noindent
Let~$D$ be the digraph which consists of a complete digraph~$A$ of order~$\lceil n/2\rceil +k/2-1$
and a complete digraph~$B$ of order~$\lfloor n/2\rfloor +k/2$ which has precisely $k-1$ vertices in
common with~$A$. It is easy to check that $\delta(D)=|A|-1=\lceil (n+k)/2\rceil -2$.
To see that~$D$ is not $k$-ordered, pick vertices $s_1,s_3,\dots,s_{k-1}$ in~$A\setminus B$
and $s_2,s_4,\dots,s_k$ in~$B\setminus A$. Then every cycle in~$D$ which
encounters $s_1,\dots,s_k$ in this order would meet~$A\cap B$ when going from~$s_i$ to~$s_{i+1}$,
i.e.~it would meet~$A\cap B$ $k$ times, which is impossible.

\medskip

\noindent
\textbf{Case 3.} \emph{$k\ge 3$ is odd and $n$ is odd.}

\smallskip

\noindent
This time we take~$D$ to be the digraph which consists of two complete digraphs~$A$ and~$B$ of
order~$(n+k)/2-1$ having~$k-2$ vertices in common. Then $\delta(D)=|A|-1=(n+k)/2 -2$.
To see that~$D$ is not $k$-ordered, pick vertices $s_1,s_3,\dots,s_k$ in~$A\setminus B$
and $s_2,s_4,\dots,s_{k-1}$ in~$B\setminus A$.
\endproof

Note that in the proof of Proposition~\ref{propbestordered} we could have omitted
the (easy) case when~$k$ is even as Proposition~\ref{propbestlinked} already gives
a digraph of the required minimum degree which is not~$k/2$-linked and thus not
$k$-ordered. 

\section{Proof of Theorem~\ref{thmkordered}}

We first prove Theorem~\ref{thmkordered} for the case%
      \COMMENT{To check the case when~$k=1$ let~$D$ be any digraph of minimum degree at least
$\lceil n/2\rceil-1$. Let~$s$ be the vertex which our cycle has to contain.
If $N^+(s)\cap N^-(s)\neq \emptyset$ then~$s$ lies on a cycle.
So suppose that $N^+(s)\cap N^-(s)= \emptyset$ and let~$A\subseteq N^+(s)$ and
$B\subseteq N^-(s)$ be sets of size $\lceil n/2\rceil-1$. Pick any vertex $a\in A$.
If~$a$ has an outneighbour~$b\in B$ then~$sab$ is a cycle. If $N^+(a)\cap B=\emptyset$
then~$n$ must be even and~$a$ must send an edge to the unique
vertex~$x$ outside $A\cup B\cup\{s\}$. Pick any vertex $b\in B$. As before, we may assume
that $N^-(b)\cap A=\emptyset$ and so~$b$ must receive an edge from~$x$.
Thus $saxb$ is a cycle.}
when~$k=2$. So suppose that~$D$ is a digraph of minimum degree at least
$\lceil n/2\rceil$. Let~$s_1$ and~$s_2$ be the vertices which our cycle has to encounter.
If $\overrightarrow{s_1s_2}$ is not an edge then $s_1,s_2\notin N^+(s_1)\cup N^-(s_2)$
and so $|N^+(s_1)\cap N^-(s_2)|\ge 2\delta(D)-(n-2)\ge 2$. Similarly, if
$\overrightarrow{s_2s_1}$ is not an edge then $|N^-(s_1)\cap N^+(s_2)|\ge 2$.
Altogether this shows that there is a cycle of length at most~4 which contains both~$s_1$
and~$s_2$.

Thus we may assume that $k \ge 3$ and that~$D$ is a digraph of minimum degree at least
$\lceil (n+k)/2\rceil-1$. Let $S:=(s_1,\dots,s_k)$ be the given sequence of vertices of~$D$
which our cycle has to encounter.
We will call these vertices \emph{special} and will sometimes also use~$S$ for
the set of these vertices. We set $s_{k+1}:=s_1$.
Given a set $I\subseteq [k]$ and a family $T:=(t_i)_{i\in I}$ of
positive integers, an \emph{$(S,I,T)$-system} is a family $(\P_i)_{i\in I}$ where
each~$\P_i$ is a set of~$t_i$ paths joining~$s_i$ to~$s_{i+1}$ and each path in~$\P_i$
has length at most~6 and is internally disjoint from~$S$, from all other paths in~$\P_i$ and from
the paths in all the other~$\P_j$.
An \emph{$(S,I)$-system} is an \emph{$(S,I,T)$-system} where~$t_i=1$ for all $i\in I$.
Thus to prove Theorem~\ref{thmkordered} we have to show that there exists an $(S,[k])$-system.

Let $I$ be the set of all those indices $i\in [k]$ for which~$D$ does not contain at
least~$6k$ internally disjoint $s_i$-$s_{i+1}$ paths of length at most~6. 

\begin{claim}\label{claim1}
It suffices to show that~$D$ contains an $(S,I)$-system.
\end{claim}

\noindent Indeed, suppose that $(\P_i)_{i\in I}$ is an $(S,I)$-system in~$D$. So each~$\P_i$ contains
precisely one path~$P_i$. We will show that for every $i\in [k]\setminus I$ we can find an
$s_i$-$s_{i+1}$ path~$P_i$ of length at most~6 which meets $S$ only in~$s_i$ and $s_{i+1}$
such that all the paths $P_1,\dots,P_k$
are internally disjoint. We will choose such a path~$P_i$ for every $i\in [k]\setminus I$
in turn. Suppose that next we want to find~$P_j$.
Recall that since $j\in [k]\setminus I$ the digraph~$D$ contains a set~$\P$ of at least~$6k$
internally disjoint $s_j$-$s_{j+1}$ paths of length at most~6.
Since at most $5(k-1)+k<6k$ vertices of $D$ lie in~$S$ or in the interior of some of the other paths~$P_i$,
one of the paths in~$\P$ must be internally disjoint from~$S$ and all the other paths~$P_i$, and so
we can take this path to be~$P_j$. This proves Claim~\ref{claim1}.

\medskip

\noindent
In order to prove the existence of an $(S,I)$-system, choose an $(S,J,T)$-system
$(\P_j)_{j\in J}$ in~$D$ such that
$J\subseteq I$ is as large as possible and subject to this $\sum_{j\in J} t_j$ is maximal.
Note that  $t_j<6k$ for all $j\in J$ since $J\subseteq I$.
Assume that $|J|<|I|$. By relabelling the special vertices, we may assume that $k\in I\setminus J$.
So we would like to extend $(\P_j)_{j\in J}$ by a suitable $s_k$-$s_1$ path.
Let $X'$ be the set of all those vertices which lie in the interior of some path belonging
to $(\P_j)_{j\in J}$. Note that
\begin{equation*}
|S\cup X'|<6k\cdot 5(k-1)+|S|<30k^2=:k_0.
\end{equation*}
Let $A:=N^+(s_k)\setminus (S\cup X')$ and $B:=N^-(s_1)\setminus (S\cup X')$.
Then
\begin{equation}\label{eqsizeAB}
|A|,|B|\ge \delta(D)-|S\cup X'|\ge n/2-k_0.
\end{equation}
Moreover, $A\cap B=\emptyset$ as otherwise
we could extend our $(S,J,T)$-system $(\P_j)_{j\in J}$ by adding the path $P_k:=s_kxs_1$
where $x\in A\cap B$, a contradiction to the choice of $(\P_j)_{j\in J}$. In particular,
this shows that the set~$X''$ of all vertices outside $A\cup B\cup S\cup X'$ has size at most~$2k_0$
and thus, setting $Y:=S\cup X'\cup X''$, we have that
$$|Y|\le 3k_0.$$
Note that~$D$
does not contain an edge $\overrightarrow{ab}$ with $a\in A$ and $b\in B$.
Indeed, otherwise we could extend $(\P_j)_{j\in J}$ by adding the
path $P_k:=s_kabs_1$.
We will often use the following claim.

\begin{claim}\label{claim2}
Let $a\in A$ and let $A'\subseteq A$ be a set of size at least~$k_0$.
Then $N^+(a)\cap A'\neq \emptyset$. Similarly, if $b\in B$ and $B'\subseteq B$ is a set of size
at least~$k_0$ then $N^-(b)\cap B'\neq \emptyset$.
\end{claim}

\noindent Suppose that $N^+(a)\cap A'= \emptyset$. Then~(\ref{eqsizeAB}) together with the
fact that~$D$ does not contain an $A$-$B$ edge implies that
$d^+(a)\le n-|B|-k_0\le n/2$, a contradiction.
The proof of the second part of the claim is similar.

\medskip

\noindent
We say that a special vertex~$s_i$ has \emph{out-type~$A$} if~$s_i$ sends at least~$k_0$ edges
to~$A$. Similarly we define when~$s_i$ has \emph{out-type~$B$}, \emph{in-type~$A$} and
\emph{in-type~$B$}. As $|Y|+2k_0\le 5k_0\le \delta(G)$,
it follows that each~$s_i$ has out-type~$A$ or out-type~$B$ (or both) and
in-type~$A$ or in-type~$B$ (or both).%
    \COMMENT{Here we use that $\delta\ge 5k_0=150k^2$, which is true if $n\ge 300k^2$
The latter is true since $k\ge 3$ and $n\ge 100k^3\ge 300k^2$.}
Note that~$s_1$ has in-type~$B$ but not in-type~$A$ whereas~$s_k$
has out-type~$A$ but not out-type~$B$. 

\begin{claim}\label{claim3}
Let $j\in J$. If~$s_j$ has out-type~$A$ then~$s_{j+1}$
has in-type~$B$ but not in-type~$A$. Similarly, if~$s_j$ has out-type~$B$ then~$s_{j+1}$
has in-type~$A$ but not in-type~$B$.
\end{claim}

\noindent Suppose that~$s_j$ has out-type~$A$ and~$s_{j+1}$ has in-type~$A$.
Let $a\in N^+_A(s_j)$. Claim~\ref{claim2} implies that $a$~sends an edge to one of the
at least~$k_0$ vertices in~$N^-_A(s_{j+1})$. Let $a'\in N^-_A(s_{j+1})$ be such
a neighbour of~$a$. Then we could extend our $(S,J,T)$-system by
adding the path $s_jaa's_{j+1}$, a contradiction.
The proof of the second part of Claim~\ref{claim3} is similar.

\begin{claim}\label{claim4}
No vertex in~$B$ sends an edge to~$A$.
\end{claim}

\noindent Suppose that $\overrightarrow{b^*a^*}$ is an edge of~$D$, where $a^*\in A$
and $b^*\in B$. Given vertices $a\in A$ and $b\in B$, put $N_{ab}:=N^+(a)\cap N^-(b)$.
Note that $N_{ab}\subseteq Y$ and $a,b\notin N^+(a)\cup N^-(b)$ as~$D$ does not
contain an $A$-$B$~edge. Thus
\begin{equation}\label{eqNab}
|N_{ab}|\ge 2\delta(D)-(n-2)=
\left(2\left\lceil\frac{n+k}{2}\right\rceil-2\right)-(n-2)\ge k.
\end{equation}
Let us now show that no special vertex~$s_i$ with $i\in J$ has out-type~$B$.
So suppose $i\in J$ and $s_i$ has out-type~$B$. Then Claim~\ref{claim3} implies that
$s_{i+1}$ has in-type~$A$. By Claim~\ref{claim2} some of the at least $k_0$ vertices
in $N^-_A(s_{i+1})$ receives an edge from~$a^*$. Let $a'$ be such a vertex.
Similarly, some of the vertices in $N^+_B(s_i)$ sends an edge to~$b^*$. Let $b'$ be such a vertex.
Then we could extend our $(S,J,T)$-system by
adding the path $s_ib'b^*a^*a's_{i+1}$, a contradiction.
This shows that whenever $s_i$ is a special
vertex of out-type~$B$  then $i\notin J$.  
Let~$Q$ denote the set of such vertices~$s_i$. 
Note that $s_k\not\in Q$ as~$s_k$ does not have out-type~$B$.
Thus each special vertex in~$Q$ forbids one index in~$J$.
 Altogether this shows that
\begin{equation}\label{eqsizeJ}
|J|\le k-1-|Q|.
\end{equation}
Let~$S_A$ be the set of all those special vertices~$s_i$ with $1\le N^-_A(s_i)<k_0$.
Let~$S_B$ be the set of all those special vertices~$s_i$ with $1\le N^+_B(s_i)<k_0$.
Let~$A^*$ be the set of all those vertices in~$A$ which do not send an edge to
some vertex in~$S_A$. Then%
     \COMMENT{$A^*$ is non-empty since $|A^*| > n/2-k_0-kk_0=n/2-30k^2-30k^3\ge
n/2-10k^3-30k^3\ge 0$ as $n\ge 80k^3$ (here we used that $k\ge 3$).}
$|A^*|> |A|-kk_0$. Similarly, let~$B^*$ be the set
of all those vertices in~$B$ which do not receive an edge from some vertex in~$S_B$. Then
$|B^*|> |B|-kk_0$.

Consider any pair~$a,b$ with $a\in A^*$ and $b\in B^*$.
As each special vertex in~$N_{ab}$ belongs to~$Q$, it follows that
\begin{equation}\label{eqNabminusS}
|N_{ab}\setminus S|\ge |N_{ab}|-|Q|\stackrel{(\ref{eqNab})}{\ge } k-|Q|
\stackrel{(\ref{eqsizeJ})}{> } |J|.
\end{equation}
Suppose first that $J\neq \emptyset$.
Given $j\in J$, let~$X''_j$ be the union of~$X''$ with the set of all vertices
lying in the interior of paths in~$\P_j$. 
As $N_{ab}\subseteq Y$ there must be an index $j_{ab}\in J$ such that~$N_{ab}$
contains at least two vertices in~$X''_{j_{ab}}$.
Note that%
     \COMMENT{$|A^*|,|B^*|> 2k_0k$ holds if $n/2-k_0-kk_0\ge 60k^3$. The latter
holds if $n/2-40k^3\ge 60k^3$, ie if $n\ge 200k^3$.}
$|A^*|,|B^*|> 2k_0k$. Thus there are $2k_0+1$ disjoint pairs~$a,b$
for which this index~$j_{ab}$ must be the same. Let $a_q,b_q$
($q=0,\dots, 2k_0$) denote these pairs and let~$j\in J$ denote the common index.

Note that $s_j$ has out-type~$A$ since we have seen before that no special vertex~$s_i$
with $i\in J$ has out-type~$B$. Claim~\ref{claim3} now implies that~$s_{j+1}$
has in-type~$B$. Pick vertices $a\in N^+_A(s_j)$ and $b\in N^-_B(s_{j+1})$
such that $a\neq a_0$ and $b\neq b_0$.
Claim~\ref{claim2} implies that there are indices $q_1,\dots,q_{k_0}$ such
that~$a$ sends an edge to each $a_{q_r}$. Apply Claim~\ref{claim2} again
to find an index~$r\le k_0$ such that~$b$ receives an edge from~$b_{q_r}$.
Let $x\in N_{a_{q_r} b_{q_r}}$ and $y\in N_{a_0b_0}$ be
distinct vertices such that $x,y\in X''_j$.
We can now modify our $(S,J,T)$-system to obtain an $(S,J\cup\{k\},T')$-system in~$D$
by replacing~$\P_j$ with the single path $s_ja a_{q_r} x b_{q_r} bs_{j+1}$
and adding the $s_k$-$s_1$ path $s_ka_0 y b_0s_1$ (see Figure~3).
If $J=\emptyset$ then we just add the $s_k$-$s_1$ path (which is still guaranteed
by~(\ref{eqNabminusS})).
In both cases this contradicts the choice
of our $(S,J,T)$-system and completes the proof of Claim~\ref{claim4}.
\begin{figure}
\centering
\psfrag{1}[][]{$A$}
\psfrag{2}[][]{$B$}
\psfrag{3}[][]{$a$}
\psfrag{4}[][]{$a_{q_r}$}
\psfrag{5}[][]{$\ \ a_0$}
\psfrag{6}[][]{$s_k$}
\psfrag{7}[][]{$\!\!\!\!s_{j+1}$}
\psfrag{8}[][]{$x$}
\psfrag{9}[][]{$y$}
\psfrag{10}[][]{$\!\!s_j$}
\psfrag{11}[][]{$\ s_1$}
\psfrag{12}[][]{$\!\!b$}
\psfrag{13}[][]{$\ \ b_{q_r}$}
\psfrag{14}[][]{$\!\!b_0$}
\includegraphics[scale=0.4]{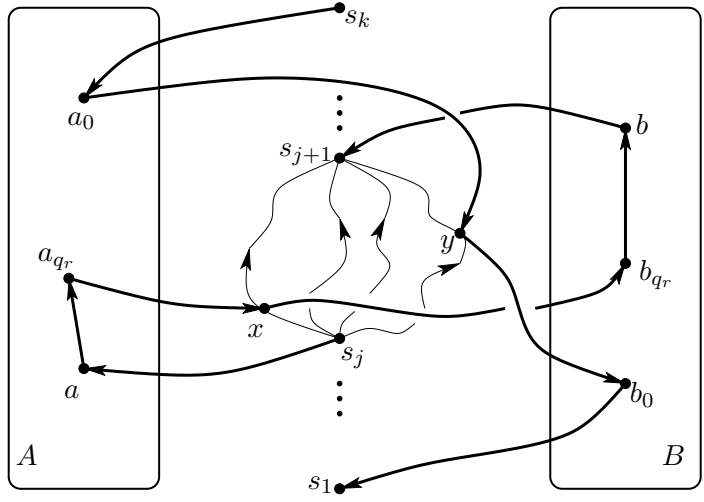}
\caption{Modifying our $(S,J,T)$-system in the proof of Claim~\ref{claim4}.}
\end{figure}%

\begin{claim}\label{claim5}
Let $a\in A$ and let $A'\subseteq A$ be a set of size at least~$k_0$.
Then $N^-(a)\cap A'\neq \emptyset$. Similarly, if $b\in B$ and $B'\subseteq B$ is a set of size
at least~$k_0$ then $N^+(b)\cap B'\neq \emptyset$.
\end{claim}

\noindent
Using Claim~\ref{claim4}, this can be shown similarly as Claim~\ref{claim2}.

\medskip

\noindent
Let $S^+_A$ be the set of all those special vertices which send an edge to~$A$
and let $S^-_A$ be the set of all those special vertices which receive an edge from~$A$.
Define~$S^+_B$ and~$S^-_B$ similarly. Note that these sets are not disjoint.
The proof of the next claim is similar to that of Claim~\ref{claim3}. (To prove the second
and third part of Claim~\ref{claim6} we use Claim~\ref{claim5} instead of Claim~\ref{claim2}.)

\begin{claim}\label{claim6}
If $j\in J$ and $s_j\in S^+_A$ then~$s_{j+1}$ cannot have in-type~$A$.
If $j-1\in J$ and $s_j\in S^-_A$ then~$s_{j-1}$ cannot have out-type~$A$.
If~$j\in J$ and $s_j\in S^+_B$ then~$s_{j+1}$ cannot have in-type~$B$.
Finally, if~$j-1\in J$ and $s_j\in S^-_B$ then~$s_{j-1}$ cannot have out-type~$B$.
\end{claim}

\noindent
Let $q^+_A:=|S^+_A|$ and define $q^-_A$, $q^+_B$ and~$q^-_B$ similarly.
Let $\Ybar:=V(D)\setminus Y=A\cup B$ and $X:=X'\cup X''=Y\setminus S$.
Consider any pair~$a,b$ with $a\in A$ and $b\in B$.
Then
\begin{equation}\label{eqN+a}
\left\lceil\frac{n+k}{2}\right\rceil-1\le |N^+(a)|\le q^-_A+|N^+_X(a)|+|N^+_{\Ybar}(a)|
\end{equation}
and
\begin{equation}\label{eqN-b}
\left\lceil\frac{n+k}{2}\right\rceil-1\le |N^-(b)|\le q^+_B+|N^-_X(b)|+|N^-_{\Ybar}(b)|.
\end{equation}
Since $N^+_{\Ybar}(a)\cap N^-_{\Ybar}(b)=\emptyset$ (as~$D$ does not contain an $A$-$B$ edge)
and $a,b\notin N^+_{\Ybar}(a)\cup N^-_{\Ybar}(b)$ we have
$$|N^+_{\Ybar}(a)|+| N^-_{\Ybar}(b)|\le |\Ybar|-2=n-|X|-k-2.
$$
Adding~(\ref{eqN+a}) and~(\ref{eqN-b}) together now gives
$$
2\left\lceil\frac{n+k}{2}\right\rceil-2\le q^-_A+q^+_B+|N^+_X(a)|+|N^-_X(b)|+n-|X|-k-2.
$$
Hence
\begin{equation}\label{eqN+aN-b}
|N^+_X(a)\cap N^-_X(b)|\ge |N^+_X(a)|+|N^-_X(b)|-|X|
\ge 2\left\lceil\frac{n+k}{2}\right\rceil-n+k-q^-_A-q^+_B 
\ge 2k-q^-_A-q^+_B.
\end{equation}
Similarly, using Claim~\ref{claim4}, one can show that
\begin{equation}\label{eqN-aN+b}
|N^-_X(a)\cap N^+_X(b)| \ge 2k-q^+_A-q^-_B.
\end{equation}
Consider any $j\in J$. Recall that by Claim~\ref{claim3} we have that either~$s_j$ has
out-type~$A$ and~$s_{j+1}$ has in-type~$B$ or~$s_j$ has out-type~$B$ and~$s_{j+1}$ has
in-type~$A$. Let $J_{AB}$ denote the set of all those indices~$j\in J$ for which the former holds
and let~$J_{BA}$ be the set of all those~$j\in J$ for which the latter holds.
Our next aim is to estimate~$j_{AB}:=|J_{AB}|$ and~$j_{BA}:=|J_{BA}|$. 
Note that Claim~\ref{claim6} implies that if $s_j\in S^+_B$ then $j\notin J_{AB}$.
As $s_k\notin S^+_B$ and $k\notin J$, this shows that
$$j_{AB}\le k-1-|S^+_B\setminus\{s_k\}|=k-1-|S^+_B|=k-1-q^+_B.
$$
Also, if $s_j\in S^-_A$ then $j-1\notin J_{AB}$ by Claim~\ref{claim6}.
As $s_1\notin S^-_A$, this shows that
$$j_{AB}\le k-1-|S^-_A\setminus\{s_1\}|=k-1-|S^-_A|=k-1-q^-_A.
$$
Adding these two inequalites gives
\begin{equation}\label{eqjAB}
j_{AB}\le k-1-\frac{q^-_A+q^+_B}{2}.
\end{equation}
In order to give an upper bound for~$j_{BA}$, note that if $s_j\in S^+_A$
then $j\notin J_{BA}$ by Claim~\ref{claim6}. Thus
$$j_{BA}\le k-1-|S^+_A\setminus\{s_k\}|\le k-q^+_A.
$$
Also, if $s_j\in S^-_B$ then $j-1\notin J_{BA}$ by Claim~\ref{claim6}.
Thus
$$j_{BA}\le k-1-|S^-_B\setminus\{s_1\}|\le k-q^-_B.
$$
Adding these two inequalites gives
\begin{equation}\label{eqjBA}
j_{BA}\le k-\frac{q^+_A+q^-_B}{2}.
\end{equation}
Our next aim is to show that~$D$ contains a $(S,J\cup\{k\})$-system. This will complete the proof
of Theorem~\ref{thmkordered} since it contradicts the choice of our $(S,J,T)$-system.
Pick distinct vertices $a_0\in A$, $a_j\in N^+_A(s_j)$ for all $j\in J_{AB}$,
$a'_j\in N^-_A(s_{j+1})$ for all $j\in J_{BA}$,
$b_0\in B$, $b_j\in N^-_B(s_{j+1})$ for all $j\in J_{AB}$ and $b'_j\in N^+_B(s_j)$ for all $j\in J_{BA}$. 
Choose%
     \COMMENT{This can be done as by~(\ref{eqN+aN-b}) $|N^+_X(a_0) \cap N^-_X(b_0)|\ge
2k-q^-_A-q^+_B\ge 2$ (the latter holds since $q^-_A\le k-1$ as $s_1\notin S^-_A$ and
$q^+_B\le k-1$ as $s_k\notin S^+_B$).} 
a vertex $x_0\in N^+_X(a_0) \cap N^-_X(b_0)$ and link~$s_k$ to~$s_1$ by the path
$Q_k:=s_ka_0x_0b_0s_1$. (This can be done since the right hand side of~(\ref{eqN+aN-b}) is at least 2.)
To find the other paths, we distinguish two cases.

\medskip

\noindent
\textbf{Case~1.} $j_{BA}\le j_{AB}$

\smallskip

\noindent
For all $j\in J_{BA}$ we pick a vertex $x_j\in N^-_X(a'_j)\cap N^+_X(b'_j)$
such that all these~$x_j$ are pairwise distinct and distinct from~$x_0$.
Inequalities~(\ref{eqN-aN+b}) and~(\ref{eqjBA}) together imply that this can be
done.%
     \COMMENT{By~(\ref{eqN-aN+b}) we need that $2k-q^+_A-q^-_B\ge j_{BA}+1$.
But (\ref{eqjBA}) implies that $2k-q^+_A-q^-_B\ge 2j_{BA}\ge j_{BA}+1$. ($j_{BA}\neq 0$ as
otherwise we don't have to do anything.)}
Inequality~(\ref{eqN+aN-b}) together with the fact that
$$2k-q^-_A-q^+_B-1-j_{BA}\stackrel{(\ref{eqjAB})}{\ge } 2j_{AB}+1-j_{BA}\ge j_{AB}+1,
$$
implies that for all $j\in J_{AB}$ we can now pick a vertex $x_j\in N^+_X(a_j)\cap N^-_X(b_j)$
such that~$x_0$ and all the $x_j$ ($j\in J$) are pairwise distinct.
If $j\in J_{AB}$ we link~$s_j$ to~$s_{j+1}$ by the path $Q_j:=s_ja_jx_jb_js_{j+1}$. 
If $j\in J_{BA}$ we link~$s_j$ to~$s_{j+1}$ by the path $Q_j:=s_jb'_jx_ja'_js_{j+1}$.
The paths $Q_j$ ($j\in J$) and~$Q_k$ are internally disjoint and have length~4,
so they form an $(S,J\cup\{k\})$-system, as required.

\medskip

\noindent
\textbf{Case~2.} $j_{BA}> j_{AB}$

\smallskip

\noindent
We proceed similiarly as in Case~1, but this time we choose the vertices
$x_j\in N^+_X(a_j)\cap N^-_X(b_j)$ for all $j\in J_{AB}$ first.
As 
$$2k-q^+_A-q^-_B-1-j_{AB}\stackrel{(\ref{eqjBA})}{\ge } 2j_{BA}-1-j_{AB}> j_{BA}-1,
$$
inequality~(\ref{eqN-aN+b}) implies that we can then pick the
vertices $x_j\in N^-_X(a'_j)\cap N^+_X(b'_j)$ for all
$j\in J_{BA}$. The paths $Q_j$ ($j\in J$) and~$Q_k$ are then defined as before.
This completes the proof of Theorem~\ref{thmkordered}.

\medskip

Note that throughout the proof, the paths we constructed always had length at most~6
(the only case where they had length exactly~6 was in the proof of Claim~\ref{claim4}).
This means that the proof can easily be translated into polynomial algorithm
so that the exponent of the running time does not depend on $k$:
We simply start with any $(S,J,T)$-system with $J \subseteq I$. Now we 
go through the steps of the proof and find a `better' $(S,J',T')$-system
with $J' \subseteq I$. 
Claim~\ref{claim1} implies that for fixed~$k$ we only need to do this a bounded 
number of times. Since the paths we need have length at most~6 and there are only a bounded
number of cases to consider in the proof, it is clear that one can find the better system in 
polynomial time with exponent independent of~$k$. Altogether this means that the problem of
finding a cycle encountering a given sequence of~$k$ vertices is fixed parameter tractable
for digraphs whose minimum degree satisfies the conditon in Theorem~\ref{thmkordered}
(where~$k$ is the fixed parameter).
The same applies to the problem of linking~$\ell$ given pairs of vertices.
In general, even the problem of deciding whether a digraph is $2$-linked is already 
NP-complete~\cite{NPlinked}. For a survey on fixed parameter tractable digraph
problems, see~\cite{Gutinsurvey}.

\section{Acknowledgement}
We are grateful to Andrew Young for reading through the manuscript.

\medskip

{\footnotesize \obeylines \parindent=0pt

Daniela K\"{u}hn, Deryk Osthus
School of Mathematics
University of Birmingham
Edgbaston
Birmingham
B15 2TT
UK
}

{\footnotesize \parindent=0pt

\it{E-mail addresses}:
\tt{\{kuehn,osthus\}@maths.bham.ac.uk}}

\end{document}